\documentclass[letterpaper, 10 pt, conference]{ieeeconf}

\usepackage{cite}

\usepackage{graphicx}
\usepackage{amssymb}
\usepackage{url}
\usepackage[cmex10]{amsmath}
\usepackage{xcolor}
\usepackage{soul}

  \usepackage[caption=false,font=footnotesize]{subfig}

\IEEEoverridecommandlockouts
\usepackage{tikz}

\renewcommand{\baselinestretch}{1}

\usepackage{multirow}

\renewcommand\Re{\operatorname{Re}}
\renewcommand\Im{\operatorname{Im}}
\renewcommand{\baselinestretch}{0.984}
\usepackage[utf8]{inputenc}

\begin{document}
%
\title{\LARGE \bf
Empirical Investigation of Non-Convexities \\ in Optimal Power Flow Problems
}


\author{Mohammad Rasoul Narimani,$^{\ast}$ Daniel K. Molzahn,$^{\dagger}$ Dan Wu,$^{\S}$ and Mariesa L. Crow$^{\ast}$%
\thanks{${\ast}$: Department of Electrical and Computer Engineering, Missouri University of Science and Technology, Rolla, MO 65401, USA, \mbox{\texttt{\{mn9t5, crow\}@mst.edu}}}%
\thanks{${\dagger}$: Energy Systems Division, Argonne National Laboratory, Argonne, IL 60439, USA, \mbox{\texttt{dmolzahn@anl.gov}}. Funding from ARPA-E, U.S. Department of Energy, award number \mbox{DE-AC02-06CH11357} is gratefully acknowledged. Views and opinions expressed herein do not necessarily state or reflect those of the United States Government or any agency thereof.}%
\thanks{${\S}$: Department of Mechanical Engineering, Massachusetts Institute of Technology, Cambridge, MA 02139, USA, \mbox{\texttt{danwumit@mit.edu}}}%
}


\maketitle

\begin{abstract}
Optimal power flow (OPF) is a central problem in the operation of electric power systems. An OPF problem optimizes a specified objective function subject to constraints imposed by both the non-linear power flow equations and engineering limits. These constraints can yield non-convex feasible spaces that result in significant computational challenges. Despite these non-convexities, local solution algorithms actually find the global optima of some practical OPF problems. This suggests that OPF problems have a range of difficulty: some problems appear to have convex or ``nearly convex'' feasible spaces in terms of the voltage magnitudes and power injections, while other problems can exhibit significant non-convexities. Understanding this range of problem difficulty is helpful for creating new test cases for algorithmic benchmarking purposes. Leveraging recently developed computational tools for exploring OPF feasible spaces, this paper first describes an empirical study that aims to characterize non-convexities for small OPF problems. This paper then proposes and analyzes several medium-size test cases that challenge a variety of solution algorithms.
\end{abstract}

\section{Introduction}
The optimal power flow (OPF) problem seeks an optimal operating point for an electric power system in terms of a specified objective function (e.g., minimizing generation cost, matching a desired voltage profile, etc.). The feasible space for an OPF problem is dictated by equality constraints corresponding to the network physics (i.e., the power flow equations) and inequality constraints determined by engineering limits on, e.g., voltage magnitudes, line flows, and generator outputs. Non-linear constraints from the power flow equations and the engineering limits can result in non-convex feasible spaces. This paper applies an empirical approach to characterize typical non-convexities that occur in OPF feasible spaces. The geometric structures characterized in this paper are based on projections of the power injections and voltage magnitudes. 

OPF problems may have multiple local optima~\cite{bukhsh_tps} and are generally NP-Hard \cite{lavaei_tps,bienstock2015nphard}, even for radial networks~\cite{NPhard}. Since first being formulated by Carpentier in 1962~\cite{carpentier}, a broad range of algorithms have been applied to solve OPF problems, including Newton-Raphson, sequential quadratic programming, interior point methods, etc.~\cite{opf_litreview1993IandII,ferc4}. Convergence of many algorithms only ensures local optimality, i.e., no feasible points in the solution's immediate neighborhood have a better objective value. Other locally optimal points may exist outside of this immediate neighborhood, some of which may have substantially better objective values.

In contrast to local solvers, global algorithms seek the lowest-cost point in the entire feasible space. Provably obtaining the global solution is relevant for many analyses, such as multi-stage and robust optimization where providing any theoretical guarantees for the overall problem requires certifying global optimality of solutions to certain subproblems~\cite{Capitanescu,pscc2018robust}. Moreover, the large scale of power systems means that even small percentage improvements in operational efficiency can have a significant aggregate impact~\cite{ferc1}, thus motivating the development of global algorithms. 

Many recently developed global algorithms employ convex relaxation techniques, which enclose the feasible space of an OPF problem in a larger convex space. Optimizing over the convex space provides a lower bound for the OPF problem's objective value, can certify OPF infeasibility, and, when the relaxation is \emph{exact}, provides the globally optimal decision variables. A variety of convex relaxations are based on semidefinite programming (SDP)~\cite{lavaei_tps,molzahn_holzer_lesieutre_demarco-large_scale_sdp_opf,molzahn_hiskens-sparse_moment_opf,josz_molzahn-complex_hierarchy} and second-order cone programming (SOCP)~\cite{coffrin2015qc,sun2015}. Recent work is surveyed in~\cite{molzahn2017survey}.

For some practical OPF problems, these convex relaxations certify that the solutions obtained by local solvers are, in fact, globally optimal (or at least very near the global optimum)~\cite{lavaei_tps,molzahn_holzer_lesieutre_demarco-large_scale_sdp_opf,molzahn_lesieutre_demarco-global_optimality_condition,molzahn_hiskens-sparse_moment_opf,coffrin2015qc,josz_molzahn-complex_hierarchy,sun2015}. There also exist challenging test cases for which local solution algorithms may fail to yield globally optimal solutions and convex relaxations have large relaxation gaps~\cite{allerton2011,bukhsh_tps,nesta}. Thus, the challenges inherent to solving OPF problems span a range of difficulties.

An OPF problem's difficulty is closely related to convexity characteristics of the problem's feasible space. The range of difficulties suggests that some OPF feasible spaces are ``nearly convex'' in terms of the voltage magnitudes and power injections, while others exhibit significant non-convexities. The development of sufficient conditions for exactness of some convex relaxation techniques~\cite{low_tutorial2} has implications for the convexity characteristics of a certain limited class of OPF problems~\cite{lavaei_geometry}. In particular, these conditions imply that portions of the feasible spaces relevant to the minimization of active power generation are convex for OPF problems that satisfy non-trivial technical conditions. Previous work also shows that the feasible spaces of a more general class of OPF problems can have significant non-convexities~\cite{hiskens2001,bernie_opfconvexity,hill2008,madani_asilomar2013,bukhsh_tps,hicss2014,polyak2017,molzahn-opf_spaces,chiang2018,molzahn-nonconvexity_search}.

Although the existing literature makes significant progress, OPF convexity characteristics are not yet fully understood. This paper leverages two recently developed computational tools to better understand non-convexities in OPF feasible spaces. The first tool is an algorithm for reliably computing discretized representations of OPF feasible spaces~\cite{molzahn-opf_spaces}. 
The second tool is a continuation algorithm that identifies multiple local optima for OPF problems~\cite{localSolutionsWu}. 

Using these tools, this paper describes an empirical analysis to better understand causes of OPF non-convexities. This analysis randomly constructs many small OPF test cases. These test cases are not directly representative of realistic power systems due to their small sizes. However, large problems may have subregions with similar features. Moreover, experience with convex relaxations of large-scale problems suggests that non-convexities are often associated with small subregions of the system \cite{molzahn_holzer_lesieutre_demarco-large_scale_sdp_opf,molzahn_hiskens-sparse_moment_opf}. Thus, exploring the characteristics of these small test cases can provide useful lessons for understanding non-convexities in large problems. After construction, the test cases are screened to identify those likely to have non-convexities using a process based on an SDP relaxation. The feasible space computation algorithm in~\cite{molzahn-opf_spaces} is applied to the screened cases to characterize their non-convexities.

Observations and test cases in~\cite{bukhsh_tps} suggest the importance of binding lower limits on voltage magnitudes and reactive power generation with regard to OPF non-convexities. All non-convexities characterized in our numerical experiment are also related to the lower limits on voltage magnitudes and reactive power generation. Our numerical experiment thus suggests that non-convexities are more frequently associated with lower limits on voltage magnitudes and reactive power generation than other constraints, at least for problems in the parameter ranges considered in our experiment.

This paper then extends the insights gained from this numerical experiment to develop challenging medium-size OPF problems based on modifications to the IEEE test cases. Modifying the system loading, voltage limits, and reactive power limits yields OPF problems where lower limits on voltage magnitudes and reactive power generation are binding. The resulting OPF problems have multiple local optima and challenge state-of-the-art convex relaxation techniques.

In addition to empirically validating the insights gained from small problems, these medium-size test cases can serve to exercise both local and global OPF solution algorithms. While convex relaxations are exact or close to exact for many previous test cases~\cite{nesta}, the medium-size test cases developed in this paper have large optimality gaps between the best-known local solutions and the bounds from the convex relaxations. In order to determine whether the optimality gaps are due to poor local optima or poor bounds, we apply both a random search technique and the continuation algorithm in~\cite{localSolutionsWu} in order to find additional local optima. This approach yields several additional local solutions and many stationary points, but none with a better objective value than that obtained via the local solver in M{\sc atpower}~\cite{matpower}. This may suggest that the optimality gaps are due to a poor bound from the relaxations, thus motivating the development of improved convex relaxation techniques.

This paper is organized as follows. Section~\ref{l:opf_overview} overviews the OPF problem. Section~\ref{l:tools} reviews computational tools for studying OPF feasible spaces. Section~\ref{l:numerical_experiment} describes the numerical experiment that is the first main contribution of this paper. Using insights from the small test cases, Section~\ref{l:ieee_modifications} presents and studies challenging OPF problems derived by modifying several IEEE test cases, which is the second main contribution of this paper. Section~\ref{l:conclusion} concludes the paper.

\section{Overview of the OPF Problem}
\label{l:opf_overview}

This section overviews the OPF problem and its SDP relaxation. Further details are provided in~\cite{ferc1,lavaei_tps,molzahn_holzer_lesieutre_demarco-large_scale_sdp_opf}. 

Consider an $n$-bus system, where $\mathcal{N} = \left\lbrace 1, \ldots, n \right\rbrace$ is the set of buses, $\mathcal{G}$ is the set of generator buses, and $\mathcal{L}$ is the set of lines. Let $\mathbf{Y}$ denote the network admittance matrix. Let $P_{Dk} + \mathbf{j} Q_{Dk}$ represent the active and reactive load demand at bus~$k\in\mathcal{N}$, where $\mathbf{j}$ is the imaginary unit. Let $V_k$ represent the voltage phasor at bus~$k\in\mathcal{N}$, with the angle of $V_1$ equal to zero to set the angle reference. Define the rank-one matrix $\mathbf{W} = V V^H \in\mathbb{H}^{n}$, where $\mathbb{H}^n$ denotes the set of $n\times n$ Hermitian matrices. Superscripts ``max'' and ``min'' denote specified upper and lower limits. Buses without generators have maximum and minimum generation set to zero. Define a convex quadratic cost of active power generation with coefficients $c_{2,k} \geq 0$, $c_{1,k}$, and $c_{0,k}$ for $k\in\mathcal{G}$. 


Each line $\left(l,m\right)\in\mathcal{L}$ is modeled by an ideal transformer with turns ratio $\tau_{lm} e^{\mathbf{j}\theta_{lm}} \colon 1$ in series with a $\Pi$ circuit with mutual admittance $y_{lm}$ and total shunt susceptance $\mathbf{j} b_{sh,lm}$. Define $e_k$ as the $k^{th}$ column of the identity matrix. Let $\overline{\left(\cdot\right)}$, $\left(\cdot\right)^\intercal$, and $\left(\cdot\right)^H$ denote the complex conjugate, transpose, and complex conjugate transpose, respectively. Define the matrices $\mathbf{H}_k = \frac{\mathbf{Y}^H e_k^{\vphantom{\intercal}} e_k^{\intercal} + e_k^{\vphantom{\intercal}} e_k^{\intercal}\mathbf{Y}}{2}$, $\tilde{\mathbf{H}}_k = \frac{\mathbf{Y}^H e_k^{\vphantom{\intercal}} e_k^{\intercal} - e_k^{\vphantom{\intercal}} e_k^{\intercal}\mathbf{Y}}{2\mathbf{j}}$, $\mathbf{F}_{lm} = \frac{1}{\tau_{lm}^2}\left(\overline{y}_{lm}-\mathbf{j}b_{sh,lm}/2\right) e_l^{\vphantom{\intercal}} e_l^{\intercal} - \overline{y}_{lm} / \left(\tau_{lm}\mathrm{e}^{-\mathbf{j}\theta_{lm}}\right) e_m^{\vphantom{\intercal}} e_l^{\intercal}$, and $\mathbf{F}_{ml} = \left(\overline{y}_{lm}-\mathbf{j}b_{sh,lm}/2\right) e_m^{\vphantom{\intercal}} e_m^{\intercal} - \overline{y}_{lm} / \left(\tau_{lm}\mathrm{e}^{\mathbf{j}\theta_{lm}}\right) e_l^{\vphantom{\intercal}} e_m^{\intercal}$.

The OPF problem is
\begin{subequations}
\small
\label{eq:opf}
\begin{align}
\nonumber
& \min_{V\in\mathbb{C}^n}\quad \sum_{k\in\mathcal{G}} c_{2,k} \left(\mathrm{tr}\left(\mathbf{H}_k \mathbf{W}\right)+P_{Dk}\right)^2 \\[-10pt] \label{eq:opf_obj}& \hspace*{60pt} + c_{1,k} \left(\mathrm{tr}\left(\mathbf{H}_k \mathbf{W}\right)+P_{Dk}\right) + c_{0,k} \hspace{-150pt}\\[-5pt] 
\nonumber & \mathrm{subject\;to} \\
\label{eq:opf_P}
& P_k^{min} \leq \mathrm{tr}\left(\mathbf{H}_k \mathbf{W}\right) + P_{Dk} \leq P_k^{max} & \forall k\in\mathcal{N}\\
\label{eq:opf_Q}
& Q_k^{min} \leq \mathrm{tr}(\tilde{\mathbf{H}}_k \mathbf{W}) + Q_{Dk} \leq Q_k^{max} & \forall k\in\mathcal{N}\\
\label{eq:opf_V2}
& \big(V_k^{min}\big)^2 \leq \mathrm{tr}\left(e_k^{\vphantom{\intercal}} e_k^{\intercal} \mathbf{W}\right) \leq \big(V_k^{max}\big)^2 & \forall k\in\mathcal{N}\\
\nonumber
\nonumber
& \!\left\lbrace\mathrm{tr}\left[ \left(\mathbf{F}_{lm}+\mathbf{F}_{lm}^H\right) \mathbf{W} \right]\right\rbrace^2 \!+\! \left\lbrace\mathrm{tr}\left[ \mathbf{j}\left(\mathbf{F}_{lm}^H-\mathbf{F}_{lm}\right) \mathbf{W} \right]\right\rbrace^2 \hspace{-55pt}\\ 
\label{eq:Sft} & \qquad \leq 4\left(S_{lm}^{max}\right)^2 & \forall \left(l,m\right) \in\mathcal{L} \\
\nonumber
& \!\left\lbrace\mathrm{tr}\left[ \left(\mathbf{F}_{ml}+\mathbf{F}_{ml}^H\right) \mathbf{W} \right]\right\rbrace^2 \!+\! \left\lbrace\mathrm{tr}\left[\mathbf{j}\left(\mathbf{F}_{ml}^H-\mathbf{F}_{ml}\right) \mathbf{W} \right]\right\rbrace^2 \hspace{-55pt}\\ 
\label{eq:Stf}
& \qquad \leq 4\left(S_{lm}^{max}\right)^2 & \forall \left(l,m\right) \mathcal{L} \\
\label{eq:opf_VVH}
& \mathbf{W} = V V^H
\end{align}
\end{subequations}
where $\mathrm{tr}\left(\cdot\right)$ is the trace. Constraints \eqref{eq:opf_P}--\eqref{eq:opf_V2} are linear in the entries of $\mathbf{W}$. The objective~\eqref{eq:opf_obj} and line flow constraints~\eqref{eq:Sft}--\eqref{eq:Stf} are convex in the entries of $\mathbf{W}$. Thus, all the non-convexity in~\eqref{eq:opf} is contained in the rank constraint~\eqref{eq:opf_VVH}.

The numerical experiment in Section~\ref{l:numerical_experiment} uses an SDP relaxation of the OPF problem as part of a screening step to identify test cases which may have relevant non-convexities. This SDP relaxation is formed by replacing~\eqref{eq:opf_VVH} with a positive semidefinite constraint $\mathbf{W} \succeq 0$~\cite{lavaei_tps}. The solution to the SDP relaxation provides a lower bound on the OPF problem's optimal objective value. If the condition $\mathrm{rank}\left(\mathbf{W}\right) = 1$ is satisfied, the lower bound provided by the SDP relaxation is \emph{exact}. Conversely, if $\mathrm{rank}\left(\mathbf{W}\right) > 1$, the lower bound may be strictly below the OPF problem's global optimum. An \emph{optimality gap} is then computed as the percent difference between the objective values for a local solution to~\eqref{eq:opf} and the lower bound from the SDP relaxation. A non-negligible optimality gap suggests the possible presence of a non-convexity in the OPF problem's feasible space near the global solution. 

Note that the OPF problem formulation~\eqref{eq:opf} does not consider some possible sources of non-convexity that are present in more general OPF problem formulations (e.g., contingency constraints, discrete devices such as switched shunts, models of uncertainty, etc.)~\cite{Capitanescu,maria2013powertech,roald2017,roald2017irep}. A variety of approaches address these possible sources of non-convexity (e.g., branch-and-bound and cutting plane methods for discrete variables~\cite{belotti2013}, chance-constrained formulations~\cite{roald2017,maria2013powertech,roald2017irep}, etc.). Many of these approaches solve the OPF formulation~\eqref{eq:opf} as a subproblem within a broader algorithm. Therefore, identifying non-convexities inherent to the OPF formulation~\eqref{eq:opf} is relevant to a wide range of problems. Future work will study the impacts of other types of OPF constraints on the feasible spaces' convexity characteristics.

\section{Tools for Studying OPF Feasible Spaces }
\label{l:tools}

This section first describes an algorithm that computes the feasible space (i.e., the set of points satisfying~\eqref{eq:opf_P}--\eqref{eq:opf_VVH}) for small OPF problems and then discusses approaches for finding multiple local optima. The numerical experiments in the following sections employ both of these algorithms to characterize OPF non-convexities.

\subsection{Computing the Feasible Spaces of Small OPF Problems}
Reference~\cite{molzahn-opf_spaces} presents an algorithm for computing a discretized representation of the feasible spaces for small OPF problems. The algorithm discretizes an OPF problem's feasible space into a set of points, each of which represents a \emph{power flow problem} (i.e., fixed voltage magnitudes at all generator buses, fixed active power injections at all generator buses except for a single ``slack'' bus which sets the angle reference, and fixed active and reactive power injections at all load buses). Observe that the expressions for power injections~\eqref{eq:opf_P}, \eqref{eq:opf_Q} and squared voltage magnitudes~\eqref{eq:opf_V2} can be written as polynomials in $V$ and $\overline{V}$ via substitution of~\eqref{eq:opf_VVH}. Expanding these complex polynomials in terms of the real and imaginary components of $V$ and $\overline{V}$ reveals a power flow formulation in terms of quadratic polynomials with real variables. See, e.g.,~\cite{molzahn-opf_spaces,mehta_molzahn_turitsyn-acc2016} for further details. 

Writing the power flow equations in a polynomial representation enables application of the ``Numerical Polynomial Homotopy Continuation'' (NPHC) algorithm, which is based on theory from algebraic geometry. The theoretical guarantees inherent to the NPHC algorithm ensure that the power flow problems corresponding to each discretization point are solved \emph{reliably}; i.e., the NPHC algorithm either returns all power flow solutions or certifies infeasibility. After solving the power flow equations corresponding to each discretization point, a screening step eliminates the solutions which fail to satisfy all the OPF problem's constraints. The remaining points are all feasible for the OPF problem, thus reliably providing a discretized representation of the entire feasible space.

The feasible space computation algorithm is only applicable to small OPF problems due to both the computational limits of the NPHC algorithm and the ``curse of dimensionality'' corresponding to the discretization of the feasible space with increasing degrees of freedom. Using convex relaxation techniques to quickly eliminate many infeasible points, the feasible space computation algorithm in \cite{molzahn-opf_spaces} is tractable for OPF problems with up to approximately ten buses and three generators. This paper's numerical experiments work within these limitations to first characterize non-convexities in small OPF problems. Lessons learned from the small problems are then applied to construct and study larger test cases. 

\subsection{Computing Multiple Local Optima}
\label{l:local_optima_search}
The presence of multiple local optima indicates the existence of non-convexities in an OPF problem's feasible space. Algorithms for computing multiple local optima therefore provide a means for investigating the associated non-convexities. Convergence of local solution algorithms depends on the selected initialization. Thus, initializing a local algorithm with various power flow solutions corresponding to random operating points is one approach for computing multiple local optima. The numerical experiments in Section~\ref{l:ieee_modifications} search for multiple local optima using at least two hundred initializations for the ``MIPS'' solver in M{\sc atpower}~\cite{matpower}. 

A more sophisticated algorithm was recently proposed in~\cite{localSolutionsWu}. Starting from a single local optimum obtained from a local solver, the algorithm in~\cite{localSolutionsWu} applies a continuation method to trace between solutions to the first-order necessary conditions for local optimality. To ensure boundedness of the continuation traces, the continuation method is applied to an ``elliptical'' representation of the first-order optimality conditions. To maintain computational tractability, we use a two-round enumeration approach; see~\cite{localSolutionsWu} for further details. This approach is capable of finding multiple local optima for problems with several tens to hundreds of buses.

\section{Investigating the Causes of OPF Non-Convexities via a Numerical Experiment}
\label{l:numerical_experiment}

The first contribution of this paper is a numerical experiment conducted to better understand the characteristics of OPF problems with non-convex feasible spaces. Specifically, this numerical experiment develops an approach for randomly constructing many small (three- to five-bus) test cases with realistic ranges for the electrical parameters. Each test case is then screened for possible non-convexities based on the \emph{optimality gap} between the objective value of a local solution and the lower bound from an SDP relaxation. The feasible spaces for the test cases identified via this screening process are then computed using that algorithm in~\cite{molzahn-opf_spaces}, which allows for characterization of the non-convexities via visual inspection. This section discusses this approach in more detail and then presents illustrative examples and various observations about the non-convexities.


\subsection{Randomly Generating and Screening Small Test Cases}
\label{l:random_procedure}

The following procedure was used to randomly construct a large number of small (three- to five-bus) OPF test cases~\cite{RTS3bus_acyclic_connected,RTS3bus_acyclic_disconnected,RTS3bus_cyclic,RTS4bus,RTS5bus}. The number of lines were sampled from a uniform distribution, with a topology developed from a random spanning tree~\cite{broder1989} augmented with additional lines whose terminal buses were randomly selected. Limits for the voltage magnitudes and angles, active and reactive power generation, load demands, line parameters, etc. were sampled from Gaussian distributions with parameters given in Appendix~\ref{a:test_case_parameters}. Test cases without sufficient generation capacity to serve the loads were discarded as trivially infeasible. See Appendix~\ref{a:test_case_parameters} for further test case construction details. A similar test case construction approach was used to study power flow problems in~\cite{acc2016}.

Computing and studying the feasible spaces for every test case is unnecessary since many of the test cases have convex or nearly convex feasible spaces that do not further this paper's goal of characterizing non-convexities. Accordingly, the following screening process was used to identify test cases which were likely to have relevant non-convexities. Using multiple random initializations, the local solver in M{\sc atpower}~\cite{matpower} was repeatedly applied to each test case. An optimality gap was then computed by comparing the lowest objective value from any initialization to the lower bound obtained from the SDP relaxation. The screening process selected test cases with large optimality gaps ($\geq 1\%$) for further analyses via the feasible space computation algorithm in~\cite{molzahn-opf_spaces}. Visualizing various projections of the feasible spaces for these test cases revealed the relevant non-convexities. The following section discusses the lessons learned from this experiment and presents instructive examples.

As a caveat for the results in the following section, note that the screening process' reliance on the lower bound from the SDP relaxation could potentially introduce bias into the selection of test cases considered for further analyses. While not observed in any related numerical experiments, it is conceptually possible that there may exist test cases with relevant non-convexities for which the SDP relaxation does not yield large optimality gaps and are therefore excluded from the empirical study. Thus, one direction for future work is to develop alternative screening processes in order to avoid any potential biases induced by the proposed approach.

\subsection{Illustrative Examples of OPF Feasible Spaces}

The empirical experiment constructed more than 10,000 test cases using the procedure in Section~\ref{l:random_procedure}, with fewer than 10 being screened for further analysis. One observation from this empirical experiment is the relatively small fraction of test cases with large optimality gaps. This suggests that relevant non-convexities (i.e., non-convexities that are near the test cases' global optima) appear to be relatively rare, at least for test cases with parameters in the ranges described in Appendix~\ref{a:test_case_parameters}. This observation aligns with previous numerical experiments indicating that the lower bound from the SDP relaxation is often close to the global optimum~\cite{molzahn2016laplacian}.  

\begin{figure}[!t]
\centering
\subfloat[One-line diagram.]{\includegraphics[width=8cm,height=8cm,keepaspectratio]{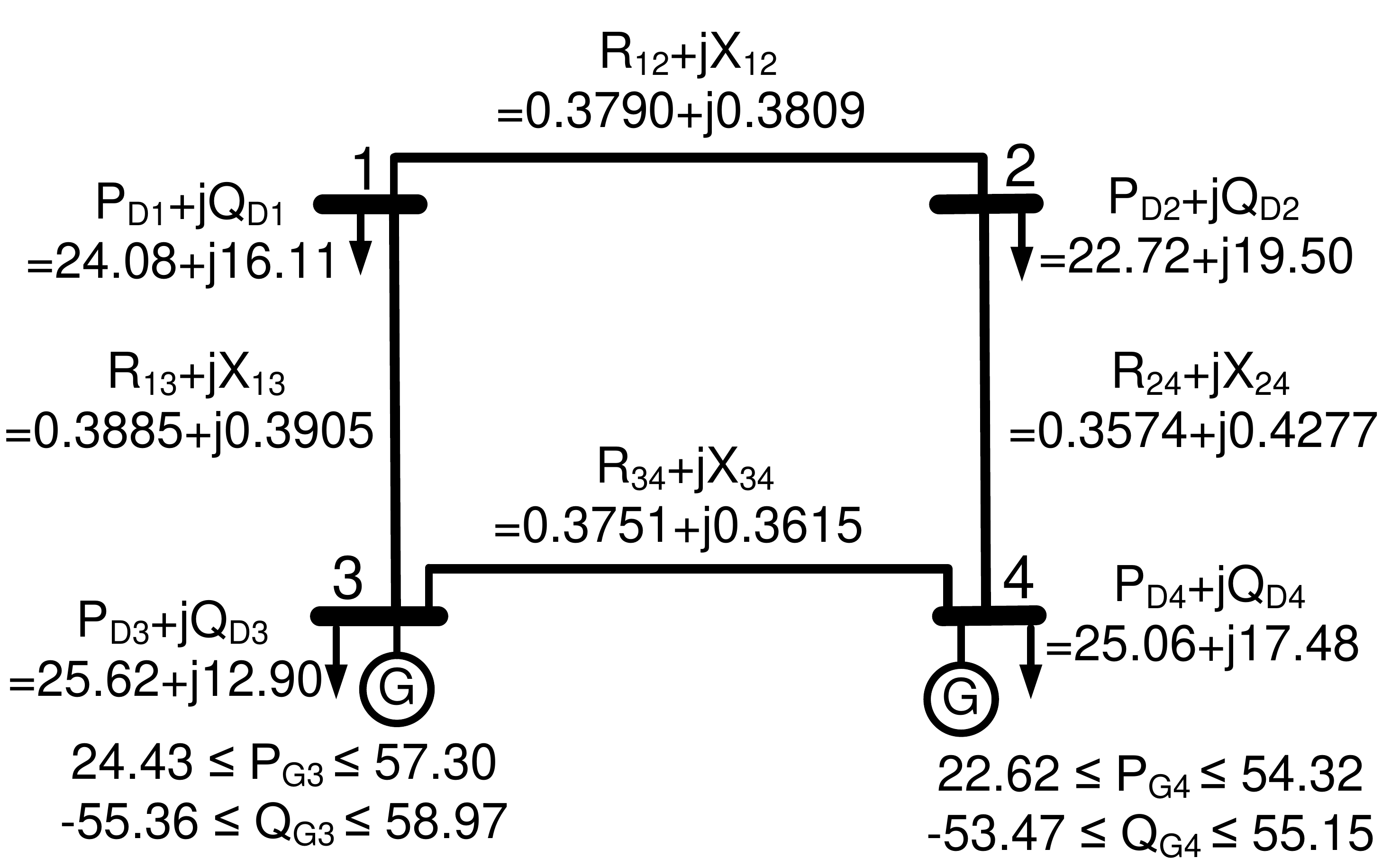}\label{f:FourBusDiagram}}\\
\subfloat[Feasible space projection.]{\includegraphics[width=8cm,height=12cm,keepaspectratio]{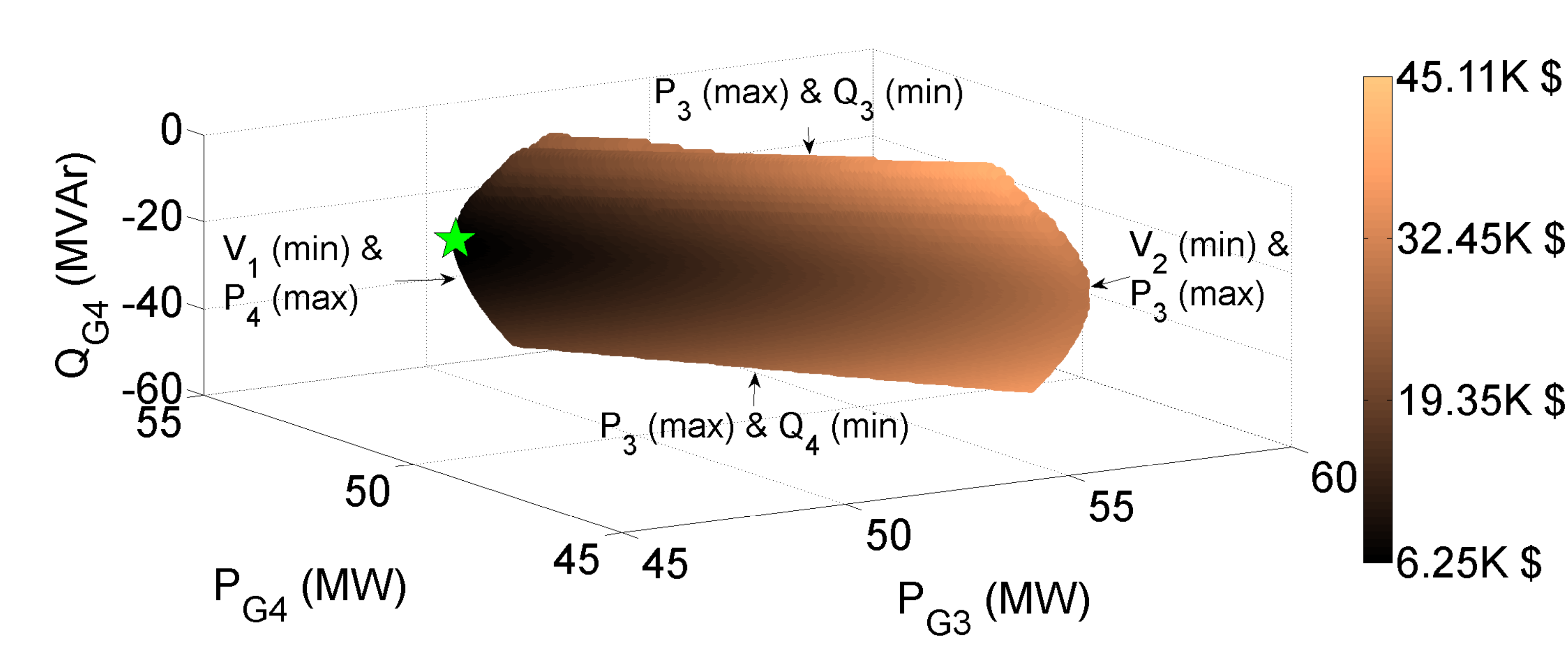}\label{f:FourBusFeasibleSpace}}
\caption{One-line diagram and feasible space projection for a ``typical'' randomly generated four-bus test case. Observe that this projection shows a convex feasible space.}
\label{f:FourBus}
\vspace*{-0.5em}
\end{figure}

\begin{figure}[!t]
\centering
\subfloat[One-line diagram.]{\includegraphics[width=9cm,height=8cm,keepaspectratio]{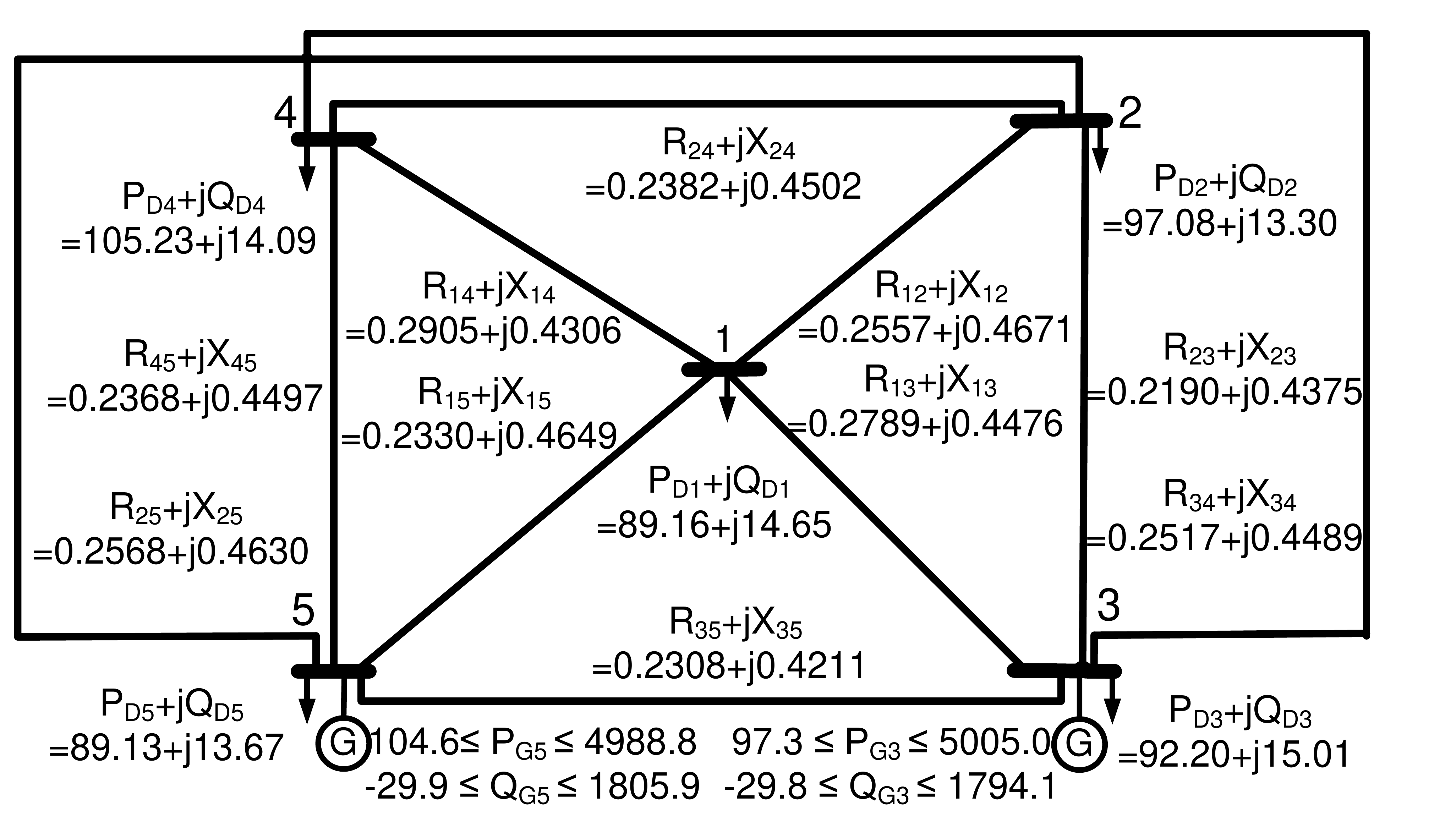}\label{f:FiveBusDiagram}}\\
\subfloat[Feasible space projection.]{\includegraphics[width=9cm,height=13cm,keepaspectratio]{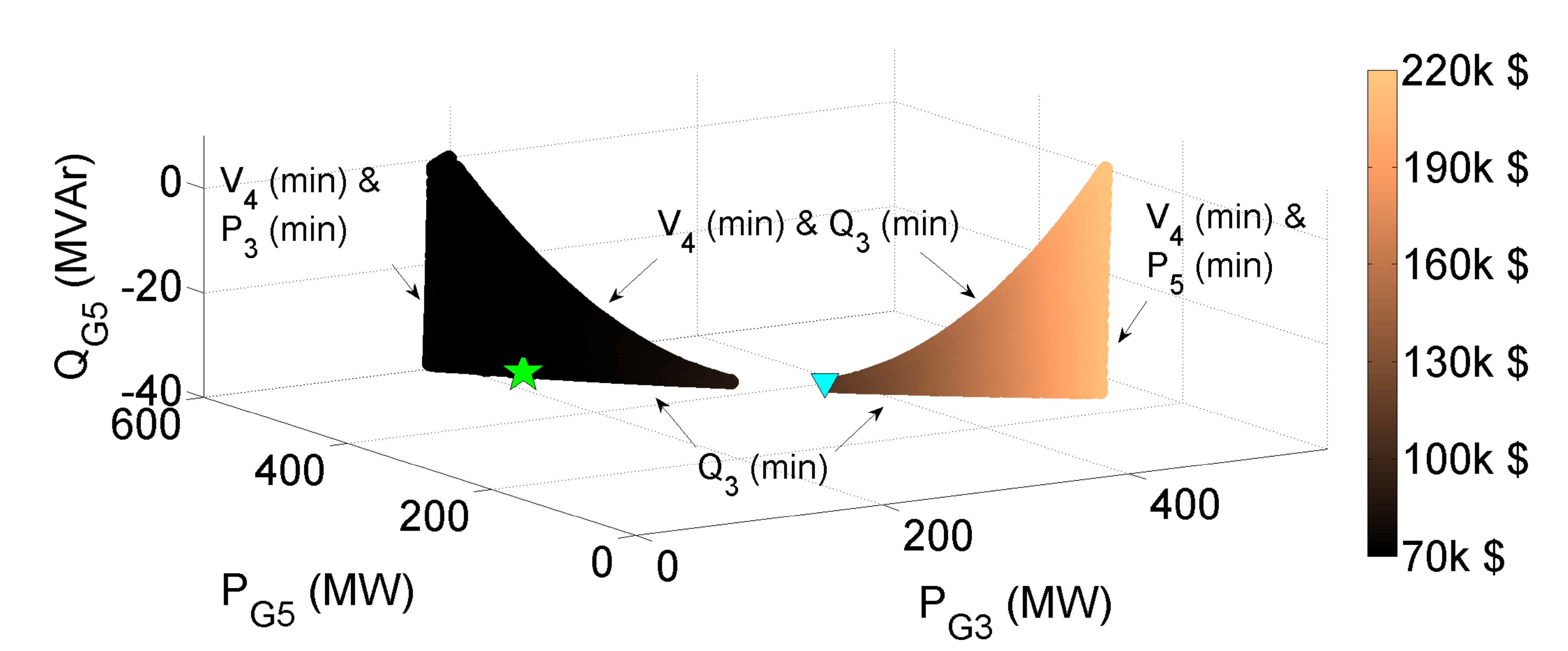}\label{f:FiveBusFeasibleSpace}}
\caption{One-line diagram and feasible space projection for a randomly generated five-bus test case. Observe that this projection shows a non-convex and disconnected feasible space.}
\label{f:FiveBus}
\vspace*{-0.5em}
\end{figure}

\begin{figure}[!t]
\centering
\vspace*{0.5em}
\subfloat[One-line diagram.]{\includegraphics[width=8.5cm,height=8.5cm,keepaspectratio]{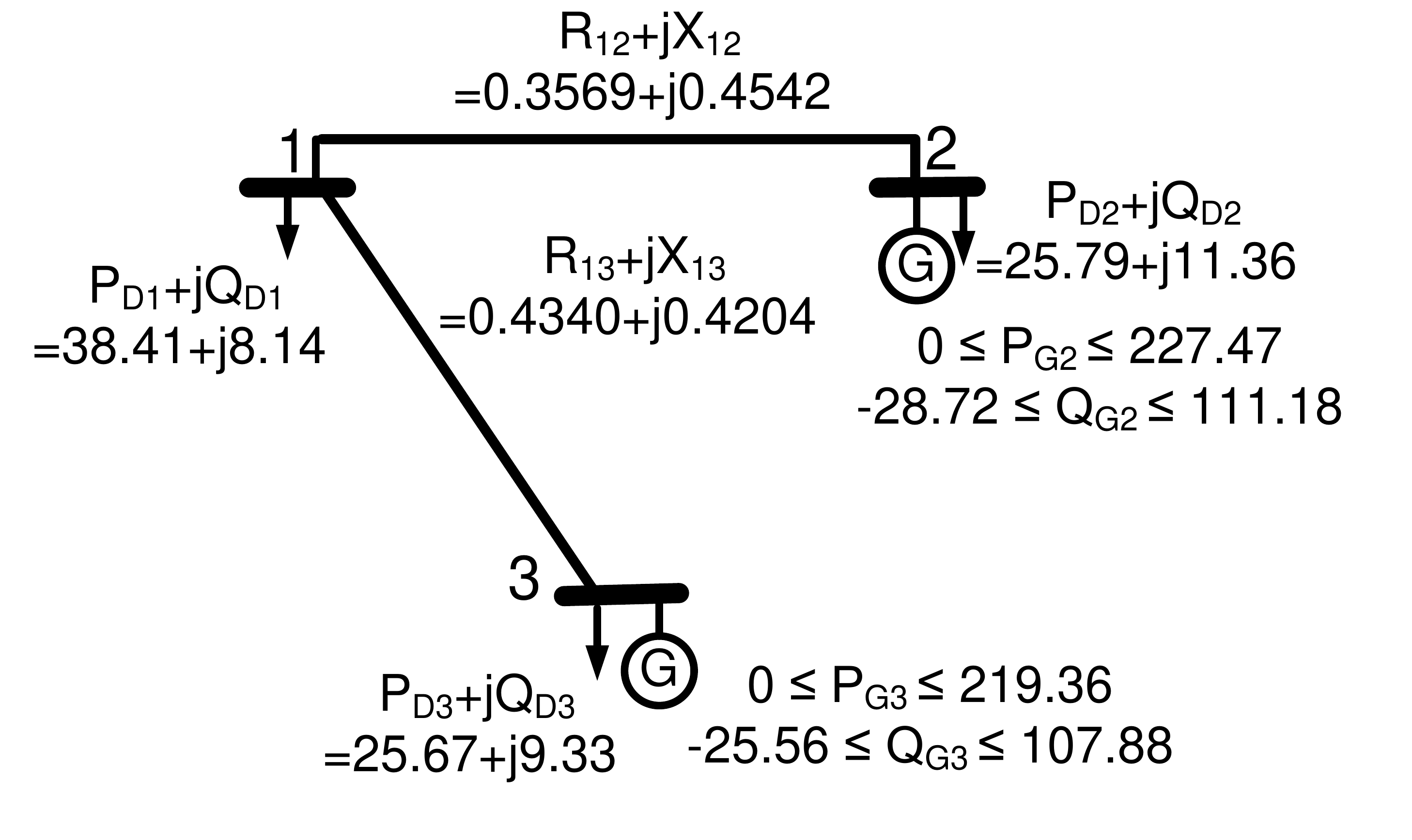}\label{f:AcyclicthreebusOnelineDiagram}}\\
\subfloat[Feasible space projection.]{\includegraphics[width=8.5cm,height=13cm,keepaspectratio]{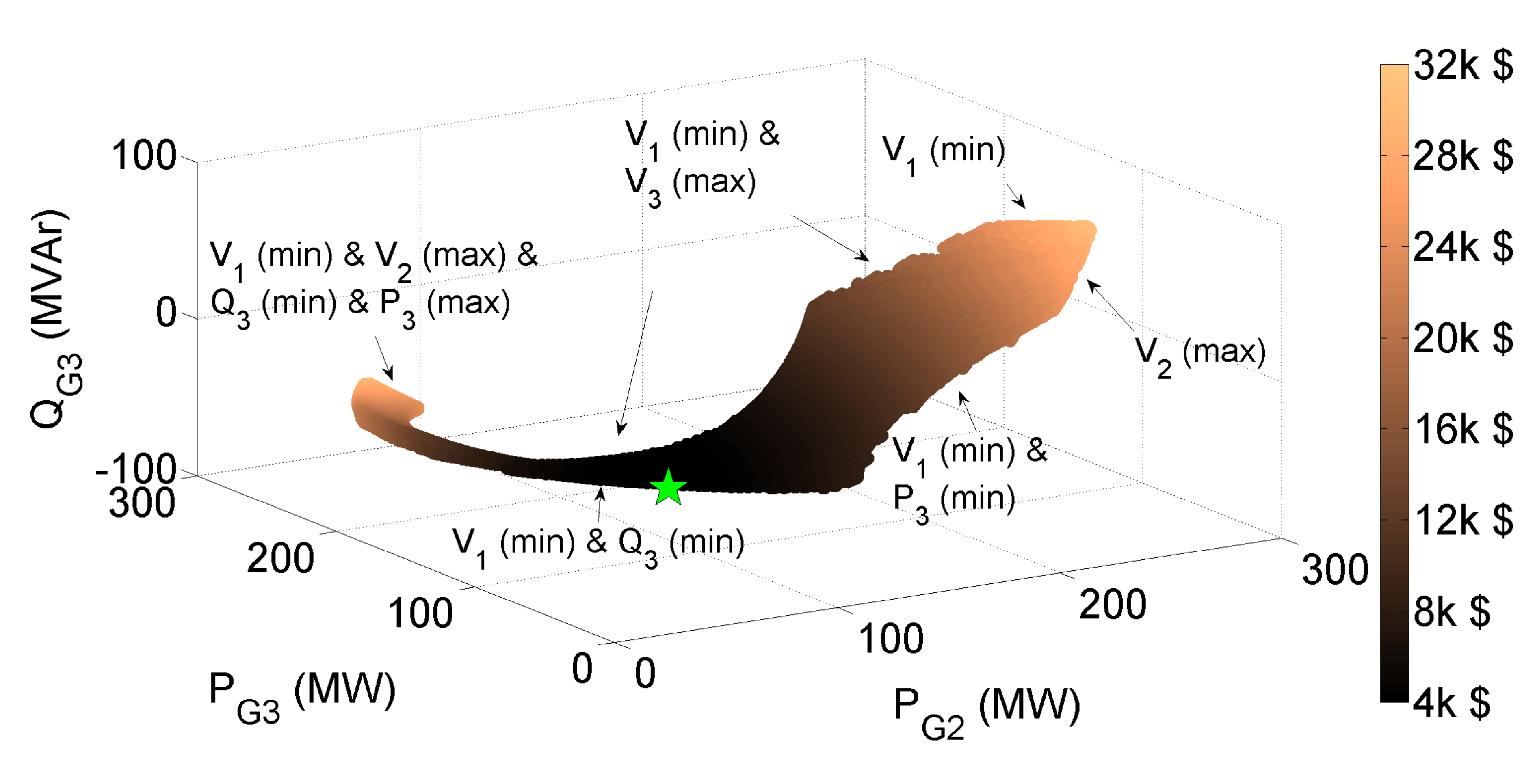}\label{f:AcyclicthreebusFS}}\\
\subfloat[Feasible space projection with tightened constraints.]{\includegraphics[width=8.5cm,height=13cm,keepaspectratio]{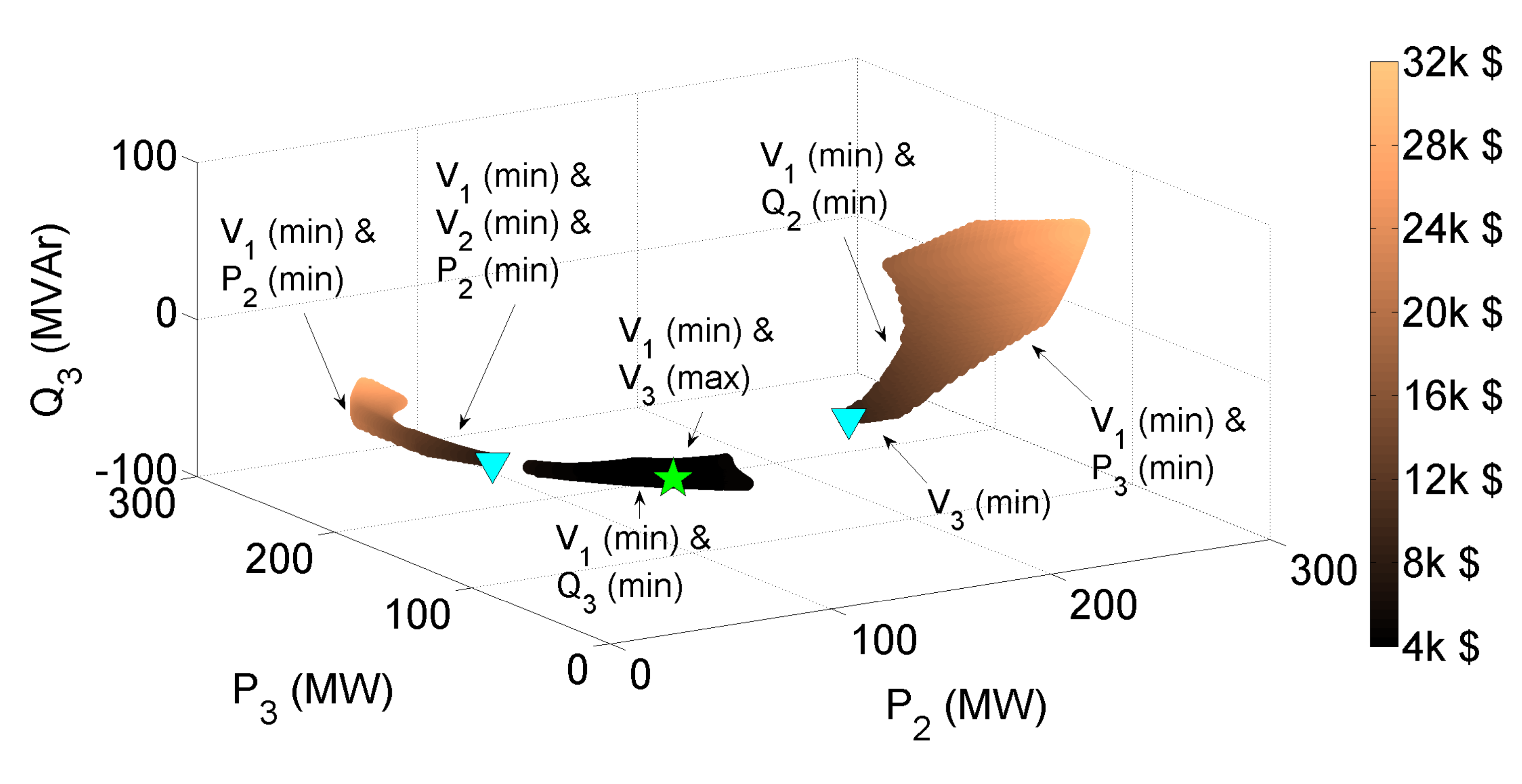}\label{f:threebus_fs_tightened}}
\caption{One-line diagram and feasible space projection for a randomly generated acyclic three-bus test case. Observe that these projections show non-convex feasible spaces. Tightening the constraints yields a disconnected feasible space in Fig.~\ref{f:threebus_fs_tightened}.}
\label{f:Acyclicthreebus}
\vspace*{-1.5em}
\end{figure}

\begin{figure}[!t]
\centering
\vspace*{0.4em}
\subfloat[One-line diagram.]{\includegraphics[width=8.5cm,height=8.8cm,keepaspectratio]{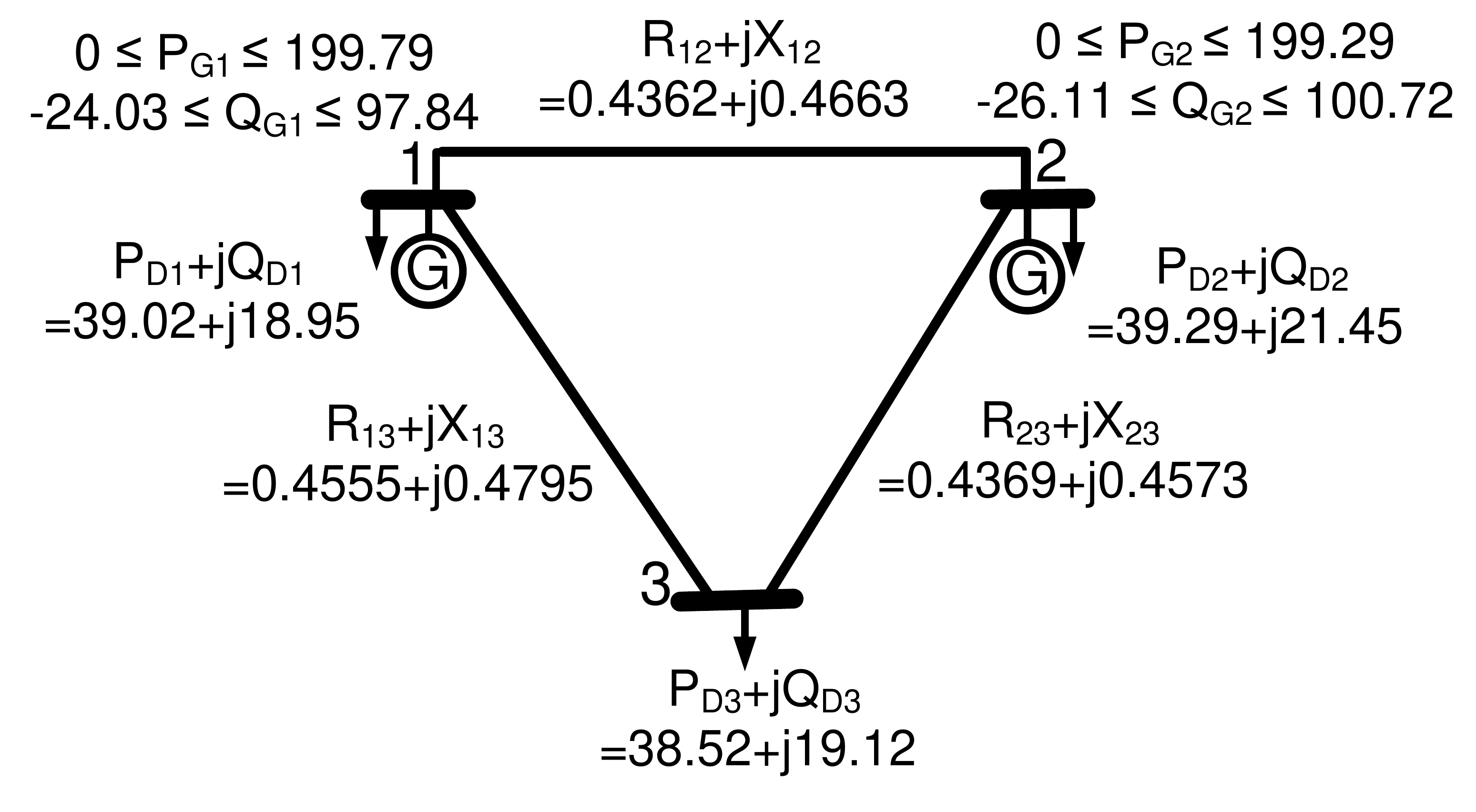}\label{f:ThreeBusMeshDiagram}}\\
\subfloat[Feasible space projection.]{\includegraphics[width=8.5cm,height=13cm,keepaspectratio]{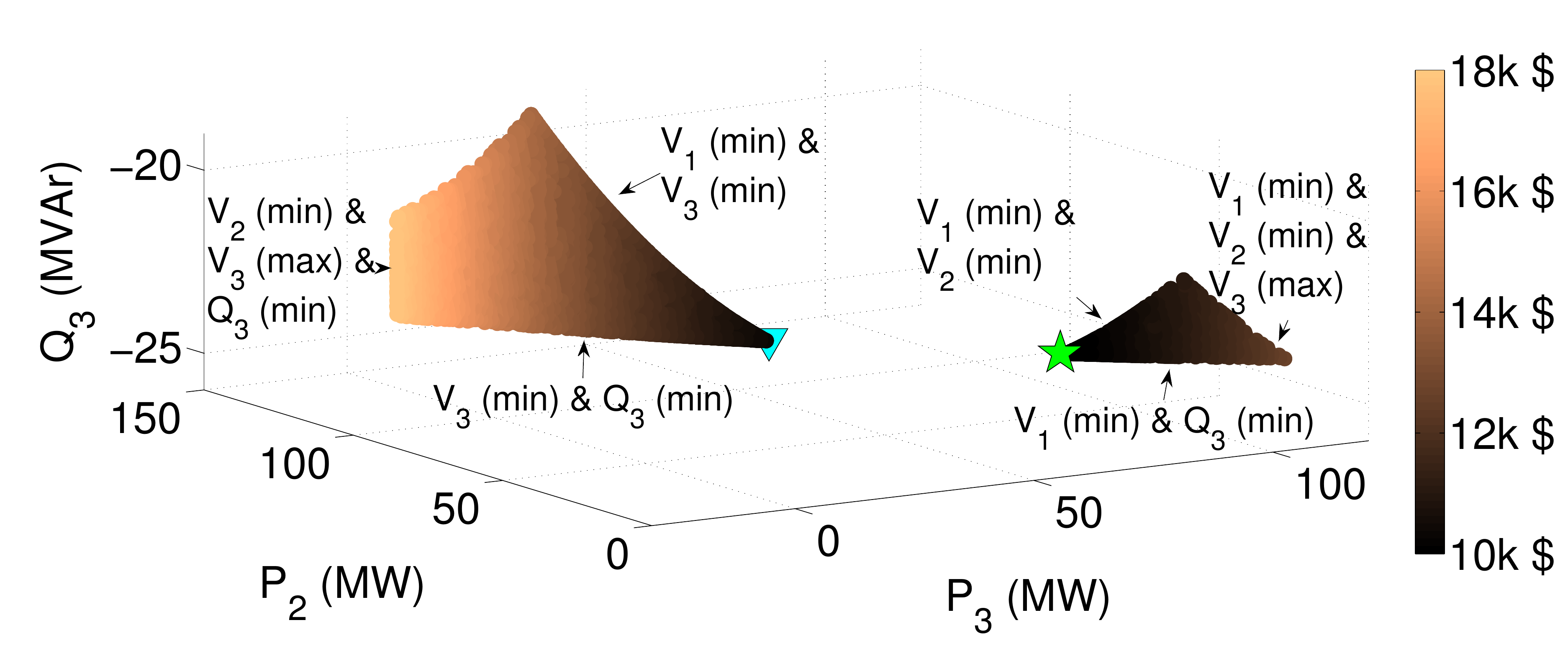}\label{f:meshed3bus_fs}}
\caption{One-line diagram and feasible space projection for a randomly generated cyclic three-bus test case. Observe that this projection shows a non-convex and disconnected feasible space.}
\label{f:MeshThreeBus}
\end{figure}

Visualizing projections of the feasible spaces for various test cases provides further insights regarding OPF non-convexities. Using the algorithm in~\cite{molzahn-opf_spaces}, this section presents several representative projections of OPF feasible spaces generated using the procedure in Section~\ref{l:random_procedure}. Figs.~\ref{f:FourBus}--\ref{f:MeshThreeBus} show one-line diagrams and projections of the corresponding feasible spaces for selected test cases. Power demands and generation ranges given in MW and MVAr. The feasible space projections are shown in terms of the active and reactive power generations (MW and MVAr) at selected buses, with the colors representing the generation cost. Line parameters are given in per unit (p.u.) on a 100~MVA base, and the shunt susceptances in the $\Pi$-circuit line model are given in Table~\ref{tab:line_shunts}. Off-nominal voltage ratios and non-zero phase shifts of transformers are tabulated in Table~\ref{tab:line_trans}. None of the flow limits are binding in any of the screened test cases. The generation cost functions and voltage magnitude limits are given in Tables~\ref{tab:gencost} and~\ref{tab:voltage_limits}, respectively.  Local and global optima are labeled with cyan triangles and green stars, respectively. 

 \begin{table}[ht]
	\caption{Line shunt values in randomly constructed test cases}
	\centering
			\begin{tabular}{ c c c c c c} \hline
 		&  \textbf{4-bus}& \textbf{5-bus}  & \textbf{3-bus (acyclic)} & \textbf{3-bus (cyclic)}\\ \hline
				$b_{1-2}$ (p.u.)& 0.3804& 0.17180 & 0.4617 & 0.4068\\ 
				$b_{1-3}$ (p.u.)& 0.4016& 0.26470 & 0.4774 & 0.4554\\ 
				$b_{2-3}$ (p.u.)& --& 0.20090  & --& 0.4376  \\
				$b_{1-4}$ (p.u.)& -- & 0.28430 & -- & --  \\ 
				$b_{1-5}$ (p.u.)& -- & 0.25632 & -- & -- \\ 
				$b_{2-4}$ (p.u.)& 0.4107& 0.02519 & -- & --\\ 
				$b_{2-5}$ (p.u.)& --& 0.21590& -- & --  \\ 
				$b_{3-4}$ (p.u.)&0.3870& 0.27260& -- & -- \\ 
				$b_{3-5}$ (p.u.)& -- & 0.20360& -- & -- \\ 
				$b_{4-5}$ (p.u.)& --& 0.28940& -- & --  \\ 
				\hline
			\end{tabular}
	\label{tab:line_shunts}
\end{table}

 \begin{table}[ht]
	\caption{Transformer details for the Five-bus test case}
	\centering
			\begin{tabular}{ c c c c } \hline
 		\textbf{Line} & \textbf{Voltage ratio} & \textbf{Phase shift (deg.)} &\\ \hline
				${1-2}$& 1.0000 & \hphantom{0}0.0000\\ 
				${1-3}$& 1.0000 & \hphantom{0}0.0000\\ 
				${3-2}$& 0.9925 & \hphantom{0}7.2099\\
				${1-4}$& 1.0000 & \hphantom{0}0.0000  \\ 
				${1-5}$& 1.0000 & \hphantom{0}0.0000 \\ 
				${2-4}$& 1.0000 & \hphantom{0}0.0000 \\ 
				${2-5}$& 1.0000 & \hphantom{0}0.0000  \\ 
				${4-3}$& 0.9950 & -2.2219 \\ 
				${3-5}$& 1.0000 & \hphantom{0}0.0000  \\ 
				${5-4}$& 1.0109 & -1.6934 \\ 
				\hline
			\end{tabular}
	\label{tab:line_trans}
\end{table}

\begin{table}[ht]
    \vspace*{1em}
	\caption{Generation cost coefficients}
	\centering
			\begin{tabular}{ c c c c} \hline
\textbf{4-bus}               &  $c_{2}$ (\$/(MW-hr)$^2$)      & $c_{1}$ (\$/(MW-hr))      & $c_{0}$ (\$)\\ \hline
\text{Generator at bus 2}    & 0.0663         & 67.2267       & 0.00\\ 
\text{Generator at bus 3}    & 0.6272         & 15.0543        & 0.00\\  \hline			
\textbf{5-bus}               &  $c_{2}$ (\$/(MW-hr)$^2$)       & $c_{1}$ (\$/(MW-hr))        & $c_{0}$ (\$)\\ \hline
\text{Generator at bus 3}    & 0.9277         & 38.7611        & 0.00\\ 
\text{Generator at bus 5}    & 0.2162         & 54.6499        & 0.40\\ \hline
\textbf{3-bus (acyclic)}     &  $c_{2}$ (\$/(MW-hr)$^2$)       & $c_{1}$ (\$/(MW-hr))        & $c_{0}$ (\$)\\ \hline
\text{Generator at bus 2}    & 0.5240         & 19.3591        & 0.00\\ 
\text{Generator at bus 3}    & 0.5480         & 16.6615        & 0.00\\  \hline
\textbf{3-bus (cyclic)}      &  $c_{2}$ (\$/(MW-hr)$^2$)      & $c_{1}$ (\$/(MW-hr))        & $c_{0}$ (\$)\\ \hline
\text{Generator at bus 2}    & 0.6408         & 49.6517        & 0.00\\ 
\text{Generator at bus 3}    & 0.6978         & 26.7824        & 0.00\\ \hline

			\end{tabular}
	\label{tab:gencost}
\end{table}

\begin{table}[ht]
	\caption{Voltage limits}
	\centering
			\begin{tabular}{ c c c } \hline
             &  $V^{max}$~(p.u.)     & $V^{min}$~(p.u.) \\ \hline
\textbf{4-bus}      & 1.10        & 0.90 \\               
\textbf{5-bus}      & 1.10         & 0.90      \\ 
\textbf{3-bus (acyclic)}    & 1.21        & 0.81 \\
\textbf{3-bus (cyclic)}     & 1.10        & 0.90 \\ 

\hline
			\end{tabular}
	\label{tab:voltage_limits}
\end{table}


Fig.~\ref{f:FourBus} shows a typical test case that did not pass the screening process (i.e., the optimality gap resulting from the SDP relaxation is small). As expected, the feasible space appears convex in terms of the power injections and voltage magnitudes. Conversely, Figs.~\ref{f:FiveBus}--\ref{f:MeshThreeBus} show examples of test cases which the screening process identified as possibly containing relevant non-convexities. The projections of the feasible spaces are indeed non-convex, with Figs.~\ref{f:FiveBusFeasibleSpace},~\ref{f:threebus_fs_tightened}, and~\ref{f:meshed3bus_fs} being disconnected. These test cases challenge a variety of optimization algorithms. Some initializations for local solvers result in convergence to suboptimal local solutions in these problems and the SDP relaxation of~\cite{lavaei_tps} is not exact.

\begin{figure}
\centering
\begin{tikzpicture}[scale=0.5]
    \draw [->,thick] (0,-4) -- (0,4) node (yaxis) [above] {$\Im\left(V_i\right)$};
    \draw [->,thick] (-4,0) -- (4,0) node (yaxis) [above right] {$\Re\left(V_i\right)$};
	\draw[even odd rule, fill=gray!30,fill opacity=0.8] circle (3) circle (2);
	\draw [->,ultra thick] (0,0) -- (2,1.5);
	
	\draw (3,0) node[below right] {\footnotesize $V_i^{max}$};
    \draw [-,thick] (3,0.25) -- (3,-0.25);
    
    \draw (2,0) node[below left] {\footnotesize $V_i^{min}$};
	\draw [-,thick] (2,0.25) -- (2,-0.25);
	
	\draw (1.8,1.5) node[above] {$V_i$};
\end{tikzpicture}
\caption{Illustration of the voltage magnitude limits~\eqref{eq:opf_V2}.}
\label{f:voltage_bounds}
\vspace*{-1.2em}
\end{figure}
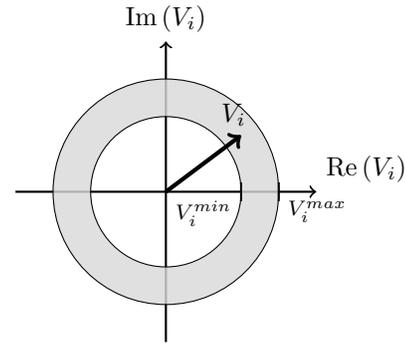

The labels in Figs.~\ref{f:FiveBusFeasibleSpace},~\ref{f:AcyclicthreebusFS},~\ref{f:threebus_fs_tightened}, and~\ref{f:meshed3bus_fs} indicate the binding limits at the boundaries of the feasible spaces. These binding limits are useful for characterizing the causes of the non-convexities. The main observation from extensive numerical experiments on these and other small test cases is that non-convexities in many OPF problems are often associated with lower limits on voltage magnitudes and reactive power generation in combination with large shunt capacitances. The lower voltage limits' relevance to non-convexities is physically intuitive: as shown in Fig.~\ref{f:voltage_bounds}, whose axes consist of the real and imaginary parts of the voltage phasors, $\Re\left(V_i\right)$ and $\Im\left(V_i\right)$, constraint~\eqref{eq:opf_V2} restricts the voltage phasors to an annulus. The lower voltage magnitude limits $\left(V_i^{min}\right)^2 \leq \Re\left(V_{i}\right)^2 + \Im\left(V_{i}\right)^2$ are thus non-convex constraints. Since increasing voltage magnitudes tends to reduce line losses, OPF problems typically have binding \emph{upper} voltage magnitude limits. In these examples, the \emph{lower} voltage limits are binding at the global optimum. To explain this, note that the large shunt capacitances in these examples result in an excess of reactive power that cannot be absorbed by the generators due to binding lower reactive power generation limits. Reducing the voltage magnitudes decreases the reactive power generated by the shunt capacitors in the lines' $\Pi$-circuit model, thus ameliorating the excess reactive power but resulting in an operating condition near the non-convexitiy associated with the lower voltage magnitude limit. 

Non-convexities were previously observed for similar operational conditions in~\cite{bukhsh_tps}. The test cases considered here thus verify the results in previous literature. Moreover, all the test cases with non-convexities characterized via the numerical experiment were associated with binding lower limits on voltage magnitudes and reactive power generation. This empirically suggests that such an operational condition is a ``common'' cause of non-convexities, at least among OPF problems with within the range of parameters considered in this experiment.

Note that the characteristics of the non-convexities (particularly disconnectedness) can be sensitive to the OPF problems' parameters. For instance, the feasible space in Fig.~\ref{f:threebus_fs_tightened} results from tightening the limits on lower reactive power generation from $Q_{G2}^{min}=-28.7$~MVAr and $Q_{G3}^{min} = -25.6$~MVAr to $Q_{G2}^{min}= -15.7$~MVAr and $Q_{G3}^{min}=-23.5$~MVAr. These modifications change this projection of the feasible space from non-convex but connected in Fig.~\ref{f:AcyclicthreebusFS} to disconnected in Fig.~\ref{f:threebus_fs_tightened}.

\section{Challenging OPF Problems Derived by Modifying IEEE Test Cases}
\label{l:ieee_modifications}

This section exploits observations from the small test cases to construct larger OPF test cases with non-convex feasible spaces. Four test cases (named ``nmwc14'', ``nmwc24'', ``nmwc57'', and ``nmwc118'' after the authors' last names and number of buses) were developed by modifying the IEEE \mbox{14-,} \mbox{24-,} \mbox{57-,} and 118-bus test cases via reducing the loading, slightly tightening the voltage limits, and significantly tightening the lower reactive power generation limits~\cite{Modified14bus,Modified24bus,Modified57bus,Modified118bus}. The goal of these modifications is to obtain test cases with operational conditions where lower limits on voltage magnitudes and reactive power generation are binding in a manner similar to the small test cases in Section~\ref{l:numerical_experiment}. Appendix~\ref{a:ieee_modifications} describes the modifications to the standard IEEE test cases provided by M{\sc atpower}~\cite{matpower}.

Applying the algorithms described in Section~\ref{l:local_optima_search} to these test cases yields multiple local optima with a wide range of objective values. Moreover, it is difficult to prove global optimality of the best known local solutions for some of these cases via relaxations with tight lower bounds, even with state-of-the-art techniques. For instance, nmwc118 has 2 local optima, and even the combination of the sparse second-order moment relaxation~\cite{molzahn_hiskens-sparse_moment_opf}, the QC relaxation~\cite{coffrin2015qc}, bound tightening~\cite{coffrin_tightening}, and a variety of related enhancements~\cite{coffrin2016strengthen_tps,coffrin2016quadtrig,sun2015} yields an optimality gap of $14.0$\%. This problem therefore appears to be particularly challenging for both traditional solvers (due to the multiple local optima) and convex relaxations. Table~\ref{tab:optimality_gap} summarizes the objective values of the known local optima and lower bounds (using a combination of the relaxations in~\cite{molzahn_hiskens-sparse_moment_opf,coffrin2015qc,coffrin_tightening,coffrin2016strengthen_tps,coffrin2016quadtrig,sun2015}) for these test cases. Note that the objective values for the local optima span wide ranges for these test cases (e.g., from \$$34664$/hr to \$$40399$/hr or equivalently from an optimality gap of $14.0$\% to $33.0$\% for nmwc118). 

 \begin{table}[h]
	\caption{Objective values for the modified IEEE test cases}
	\centering
		\begin{tabular}{ l c c c c c } \hline
		    \multicolumn{1}{c}{\textbf{Case}} & \multicolumn{2}{c}{\textbf{Local Optima (\$/hr)}} & \textbf{Lower Bound} & \textbf{Optimality} \\ 
		    \multicolumn{1}{c}{\textbf{Name}} & \textbf{Worst} & \textbf{Best} & \textbf{(\$/hr)} & \textbf{Gap} \\\hline
			 nmwc14  & \hphantom{0}3024.46& \hphantom{0}2529.87 & \hphantom{0}2529.49 & \hphantom{0}0.01\%\\
			 nmwc24  & 42667.26& 39773.02&39773.02& \hphantom{0}0.00\%\\  
			 nmwc57  & \hphantom{0}9186.12& \hphantom{0}9128.72 & \hphantom{0}9030.70 & \hphantom{0}1.09\%\\
			 nmwc118 & 40399.17& 34663.69& 30413.10 &14.00\%\\
			\hline
		\end{tabular}
	\label{tab:optimality_gap}
\end{table}

\section{Conclusion}
\label{l:conclusion}

Despite significant recent progress, there remain problems whose non-convexities challenge state-of-the-art OPF solution algorithms. Better understanding these non-convexities is important for further improving solution algorithms as well as for developing additional challenging test cases. The numerical experiment described in this paper provides a key observation regarding OPF non-convexities: all of the non-convexities identified in the numerical experiment are associated with binding lower bounds on voltage magnitudes and reactive power generation. Exploiting this observation, this paper proposes several new test cases derived by modifying the loading and tightening certain constraints in the IEEE test cases. With many local optima and large optimality gaps, these cases challenge a variety of state-of-the-art algorithms. 

Ongoing work aims to construct larger test cases which exhibit a range of difficulties. Other ongoing work involves studying the feasible spaces for test cases from the NESTA archive~\cite{nesta} with large optimality gaps. Related future work will exploit the observations in this paper to improve convex relaxation algorithms.

\appendices
\section{Parameters for Random Test Case Generation Algorithm}
\label{a:test_case_parameters}

Table~\ref{tab:table_one} provides the parameters used in the empirical experiment in Section~\ref{l:numerical_experiment}. For each test case, impedance $R+\mathbf{j}X$ (yielding admittance $g+\mathbf{j}b = 1/\left(R+\mathbf{j}X\right)$) and shunt susceptance $b$ values for the lines' $\Pi$-model equivalent circuits were randomly sampled from Gaussian distributions with mean and standard deviation of $\mu_R$, $\sigma_R$; $\mu_X$, $\sigma_X$; and $\mu_b$, $\sigma_b$, respectively, in per unit using a $100$~MVA base, with any negative values sampled for line resistances instead set to zero. Lines had an $8$\% probability of being transformers with tap ratio $\tau$ and phase-shift $\theta$ sampled from a Gaussian distribution with mean and standard deviation values of $\mu_{\tau}$, $\sigma_{\tau}$  per unit and $\mu_{\theta}$, $\sigma_{\theta}$, respectively. A bus was specified to be a generator with $30$\% probability, with the first generator selected as the reference bus. If no buses were selected to be generators, a random bus was assigned a generator and chosen to provide the angle reference. The active power injections were sampled from Gaussian distribution with mean and standard deviation values of $\mu_{P_g}$ and $\sigma_{P_g}$. Loads have a constant active and reactive power component sampled from a Gaussian distribution with mean and standard deviation $\mu_{P_d}$, $\sigma_{P_d}$ and $\mu_{Q_d}$, $\sigma_{Q_d}$, respectively. A variety of numerical experiments not detailed in this paper tested different ranges of parameter values. The parameters in Table~\ref{tab:table_one} were chosen such that the resulting test cases tend to be feasible and provide at least some examples which passed the screening process discussed in Section~\ref{l:random_procedure}.

 \begin{table}[t]
    \vspace*{1em}
	\caption{Means and standard deviations for parameter values in the randomly constructed test cases}
	\centering
		\resizebox{\columnwidth}{!}{%
			\begin{tabular}{ l c c c c c} \hline
			 &\textbf{4-bus}	& \textbf{5-bus} & \textbf{3-bus (acyclic)} & \textbf{3-bus (cyclic)}   \\ \hline
				$\mu_R$ (p.u.)& \hphantom{00}0.37& \hphantom{000}0.25& \hphantom{00}0.40& \hphantom{00}0.43 \\ 
				$\sigma_R$ (p.u.)& \hphantom{00}0.02& \hphantom{000}0.01& \hphantom{00}0.05& \hphantom{00}0.02 \\ 
				$\mu_X$ (p.u.)& \hphantom{00}0.38& \hphantom{000}0.44& \hphantom{00}0.44 & \hphantom{00}0.46\\ 
				$\sigma_X$ (p.u.)& \hphantom{00}0.02 & \hphantom{000}0.01& \hphantom{00}0.01& \hphantom{00}0.01 \\ 
				$\mu_b$ (p.u.)& \hphantom{00}0.38& \hphantom{000}0.22& \hphantom{00}0.45& \hphantom{00}0.43 \\ 
				$\sigma_b$ (p.u.)& \hphantom{00}0.02  & \hphantom{000}0.02& \hphantom{00}0.01& \hphantom{00}0.01 \\ 
				$\mu_{\tau}$ (p.u.)& \hphantom{00}0.00  & \hphantom{000}1.00& \hphantom{00}1.00& \hphantom{00}1.00  \\ 
				$\sigma_{\tau}$ (p.u.)& \hphantom{00}0.00   & \hphantom{000}0.01& \hphantom{00}0.00& \hphantom{00}0.00  \\ 
				$\mu_{\theta}$ (deg)& \hphantom{00}0.00  & \hphantom{000}0.00& \hphantom{00}0.00& \hphantom{00}0.00  \\ 
				$\sigma_{\theta}$ (deg.)& \hphantom{00}0.00  & \hphantom{000}3.00& \hphantom{00}0.00& \hphantom{00}0.00 \\ 
				$\mu_{P_{g,max}}$ (MW)& \hphantom{0}24.00 & 5000.00& 220.00& 200.00  \\ 
				$\sigma_{P_{g,max}}$ (MW)& \hphantom{00}1.00 & \hphantom{000}5.00& \hphantom{00}2.00& \hphantom{00}1.00 \\ 
				$\mu_{P_{g,min}}$ (MW)& \hphantom{0}23.00 & \hphantom{0}100.00& \hphantom{00}0.00& \hphantom{00}0.00  \\ 
				$\sigma_{P_{g,min}}$ (MW)& \hphantom{00}1.00 & \hphantom{000}2.00& \hphantom{00}0.00& \hphantom{00}0.00  \\ 
				$\mu_{Q_{g,max}}$ (MVAr)& \hphantom{0}57.00 & 1800.00& 110.00& 100.00 \\ 
				$\sigma_{Q_{g,max}}$ (MVAr)& \hphantom{00}2.00 & \hphantom{000}5.00& \hphantom{00}2.00& \hphantom{00}2.00  \\ 
				$\mu_{Q_{g,min}}$ (MVAr)& -54.00 & \hphantom{0}-30.00& -26.00& -25.00  \\ 
				$\sigma_{Q_{g,min}}$ (MVAr)& \hphantom{00}1.00  & \hphantom{000}1.00& \hphantom{00}1.00& \hphantom{00}1.00 \\ 
				$\mu_{P_d}$ (MW) & \hphantom{0}23.00 & \hphantom{00}95.00& \hphantom{0}30.00& \hphantom{0}39.00 \\ 
				$\sigma_{P_d}$ (MW)& \hphantom{00}3.00  & \hphantom{000}5.00& \hphantom{00}5.00& \hphantom{00}0.00  \\ 
				$\mu_{Q_d}$ (MVAr)& \hphantom{0}16.00  & \hphantom{00}14.00& \hphantom{0}10.00& \hphantom{0}20.00 \\ 
				$\sigma_{Q_d}$ (MVAr)& \hphantom{00}3.00  & \hphantom{000}1.00& \hphantom{00}1.00& \hphantom{00}1.00 \\ 
				\hline
			\end{tabular}
		}%
	\label{tab:table_one}
    \vspace*{-0.25em}
\end{table}

\section{Description of the Modified IEEE Test Cases}
\label{a:ieee_modifications}
Table~\ref{tab:table_two} provides the percentage changes applied to each of the IEEE test cases. Modifications to the IEEE test cases consist of decreasing active and reactive loads by $\delta_{P_d}$ and $\delta_{Q_d}$, tightening upper and lower bounds on voltage by $\delta_{\overline{V}}$ and $\delta_{\underline{V}}$, and tightening the lower bound on reactive power by $\delta_{Q_G}$.
 
 \begin{table}[t]
    \vspace*{1em}
 	\caption{Descriptions of modifications to the IEEE test systems}
 	\centering
 			\begin{tabular}{ l c c c c c} \hline
 				 & \textbf{14-bus} & \textbf{24-bus} & \textbf{57-bus} & \textbf{118-bus} \\ \hline
 				$\delta_{P_d}$~(\%) & 60.00& 55.00& 72.00& 71.00\\ 
 				$\delta_{Q_d}$~(\%) & 60.00& 55.00& 72.00& 71.00\\ 
 				$\delta_{\overline{V}}$~(\%) & \hphantom{0}0.06& \hphantom{0}0.73& \hphantom{0}0.06& \hphantom{0}0.06\\ 
 				$\delta_{\underline{V}}$~(\%) & \hphantom{0}0.06& \hphantom{0}0.73& \hphantom{0}0.06& \hphantom{0}0.06 \\ 
 				$\delta_{Q_G}$~(\%) & 95.00& 90.00& 95.00& 95.00\\ \hline
 			\end{tabular}
 	\label{tab:table_two}
    \vspace*{-1.2em}
 \end{table}





\bibliographystyle{IEEEtran}
\bibliography{IEEEabrv,ref}

\end{document}